\documentclass[12pt]{article}
\usepackage{latexsym,amsfonts,amssymb}
\usepackage[dvips]{color}
\usepackage{colordvi,multicol}

\newtheorem{thm}{Theorem}[section]
\newtheorem{cor}[thm]{Corollary}

\newtheorem{pro}[thm]{Proposition}
\newtheorem{defi}[thm]{Definition}
\newtheorem{rem}[thm]{Remark}
\newtheorem{exa}[thm]{Example}

\begin{document}
\title{\bf  Toeplitz Lemma, Complete Convergence and Complete Moment Convergence}
\author{Jiyanglin Li and Ze-Chun Hu\thanks{Corresponding
author: Department of Mathematics, Nanjing University, Nanjing
210093,  China\vskip 0cm E-mail address: huzc@nju.edu.cn}\\
 {\small Nanjing University, Nanjing, China}}

\maketitle
\date{}

\begin{abstract}
In this paper, we study the Toeplitz lemma, the  Ces\`{a}ro mean convergence theorem and the Kronecker lemma. At first,
we study ``complete convergence" versions of the Toeplitz lemma, the
Ces\`{a}ro mean convergence theorem and the Kronecker lemma. Two counterexamples show that
they can fail in general and some sufficient conditions for ``complete convergence" version of the Ces\`{a}ro mean convergence theorem are given. Secondly,
  we introduce two classes of complete moment convergence, which are stronger versions of mean convergence and consider the Toeplitz lemma,
the Ces\`{a}ro mean convergence theorem, and the Kronecker lemma
under these two classes of complete moment convergence.

\end{abstract}



\noindent

\section{Introduction}
The Toeplitz lemma and its two corollaries (the Ces\`{a}ro mean convergence theorem and the Kronecker lemma) are useful tools in the study of
probability limit theorems. For the reader's convenience, we spell out them in the following and their proofs may be found in Lo\`{e}ve (1977).

\begin{thm} {\rm (Toeplitz lemma)} Let $\{a_{nk},1\leq k\leq k_n,n\ge1 \}$ be a double array of real numbers such that for any $k\geq 1, \lim_{n\to\infty}a_{nk}=0$ and
$\sup_{n\geq 1}\sum_{k=1}^{k_n}|a_{nk}|<\infty$. Let $\{x_n,n\geq 1\}$ ba a sequence of real numbers.\\
(i) If $\lim_{n\to\infty}x_n=0$, then $\lim_{n\to\infty}\sum_{k=1}^{k_n}a_{nk}x_k=0.$\\
(ii) If $\lim_{n\to\infty}x_n=x\in \mathbf{R}$  and $\lim_{n\to\infty}\sum_{k=1}^{k_n}a_{nk}=1$, then $\lim_{n\to\infty}\sum_{k=1}^{k_n}a_{nk}x_k=x$.
\end{thm}

\begin{cor}{\rm (Ces\`{a}ro mean convergence theorem)} Let $\{x_n,n\geq 1\}$ be a sequence of real numbers and let $\bar{x}_n=\sum_{k=1}^nx_k/n,n\geq 1$. If $\lim_{n\to\infty}x_n=x\in \mathbf{R}$, then
$\lim_{n\to\infty}\bar{x}_n=x$.
\end{cor}

\begin{cor}{\rm (Kronecker lemma)} Let $\{x_n,n\geq 1\}$ and $\{b_n,n\geq 1\}$ be sequences of real numbers such that $0<b_n\uparrow\infty$. If the series $\sum_{k=1}^{\infty}x_k/b_k$ converges, then
$\lim_{n\to\infty}\frac{1}{b_n}\sum_{k=1}^nx_k=0$.
\end{cor}

By the definition of almost sure (a.s.) convergence, we know that the Toeplitz lemma and its two corollaries (the Ces\`{a}ro mean convergence theorem and the Kronecker lemma)  still hold when the numerical sequence $\{x_n,n\geq 1\}$ and real number $x$ are replaced by a sequence of random variable $\{X_n,n\geq 1\}$ and a random variable $X$, respectively, and the limit is taken to be a.s. convergence.

Recently,    Linero  and Rosalsky (2013) showed among other things that ``convergence in probability" versions of the Toeplitz lemma, the Ces\`{a}ro mean convergence theorem and the  Kronecker lemma can fail, and their ``mean convergence" versions are true.

We know that both a.s. convergence and mean convergence imply convergence in probability. Then there are the following two questions:

{\it Question 1. Do the Toeplitz lemma, the Ces\`{a}ro mean convergence theorem and the  Kronecker lemma hold under stronger convergence than a.s. convergence?}

{\it Question 2. Do the Toeplitz lemma, the Ces\`{a}ro mean convergence theorem and the  Kronecker lemma hold under stronger convergence than mean convergence?}

In Section 2, we study Question 1 and consider ``complete convergence " versions of the Toeplitz lemma, the Ces\`{a}ro mean convergence theorem and the Kronecker lemma. At first, we will give two examples to show that
they can fail in general. Then we give some sufficient conditions for the Ces\`{a}ro mean convergence theorem under complete convergence.

Let $\{X, X_n,n\geq 1\}$  be a sequence of random variables on some probability space $(\Omega,\mathcal{F},P)$. If  $\forall \varepsilon>0$,
$$
\sum_{n=1}^{\infty}P\{|X_n-X|\geq \varepsilon\}<\infty,
$$
then  $\{X_n,n\geq 1\}$  is said to converge completely to $X$ (write $X_n\stackrel{c.c.}{\longrightarrow}X$, or $X_n\to X$ c.c. for short).  This concept was introduced by Hsu and Robbins (1947).
Let $\{X,X_n,n\geq 1\}$ be a sequence of independent and indentically distributed (i.i.d.) random variables and set $S_n=\sum_{k=1}^nX_k,n\geq 1$. Hsu and Robbins (1947) proved that if $E[X]=0$ and $E[X^2]<\infty$, then $S_n/n\stackrel{c.c.}{\longrightarrow}0$. The converse was proved by Erd\"{o}s (1949, 1950). The Hsu-Robbins-Erd\"{o}s theorem was generalized in various ways, see, Baum and Katz (1965), Gut  (1978, 1980), Li et al. (1995), Lanzinger (1998), Sung and Volodin (2006), Sung (2007), Gut and Stadtm\"{u}ller (2011), and Chen and Sung (2014).

In Section 3, we study Question 2. To that end, in view of the relations between convergence in probability and complete convergence, we introduce two classes of complete moment convergences, which are stronger versions of mean convergence.
Let $p>0$.

\begin{defi}
$\{X_n,n\geq 1\}$  is said to s-$L^p$ converge to $X$ (denote $X_n\stackrel{s\mbox{-}L^p}{\longrightarrow} X$  for short), if
\begin{eqnarray*}
\sum_{n=1}^{\infty}E[|X_n-X|^p]<\infty.
\end{eqnarray*}
\end{defi}

\begin{defi}
$\{X_n,n\geq 1\}$  is said to s$^*$-$L^p$ converge to $X$ (denote $X_n\stackrel{s^*\mbox{-}L^p}{\longrightarrow} X$  for short), if
$$
\sum_{n=1}^{\infty}\|X_n-X\|_p<\infty,
$$
\end{defi}
where $\|X_n-X\|_p=(E[|X_n-X|^p])^{1/p}$.

\begin{rem}
(i) Obviously, if $X_n\stackrel{s\mbox{-}L^p}{\longrightarrow} X$ or $X_n\stackrel{s^*\mbox{-}L^p}{\longrightarrow} X$ for some $p>0$, then $\|X_n-X\|_p\to 0$.\\
(ii) By Markov's inequality, we know that if $X_n\stackrel{s\mbox{-}L^p}{\longrightarrow} X$ for some $p>0$, then $X_n\stackrel{c.c.}{\longrightarrow}X$ and thus $X_n\stackrel{a.s.}{\longrightarrow}X$ by the Borel-Cantelli lemma.\\
(iii) If $p>1$ and $X_n\stackrel{s^*\mbox{-}L^p}{\longrightarrow} X$, then $X_n\stackrel{s\mbox{-}L^p}{\longrightarrow} X$; if $0<p<1$ and $X_n\stackrel{s\mbox{-}L^p}{\longrightarrow} X$, then $X_n\stackrel{s^*\mbox{-}L^p}{\longrightarrow} X$.
\end{rem}

Chow (1988) first investigated the complete moment convergence, and obtained the following result. Let $\{X,X_n,n\geq 1\}$ be a sequence of i.i.d. random variables with $E[X]=0$. Let $1\leq p<2$ and $\gamma\geq p$. If $E[|X|^{\gamma}+|X|\log(1+|X|)]<\infty$, then
\begin{eqnarray}\label{Chow-a}
\sum_{n\geq 1}n^{\frac{\gamma}{p}-2-\frac{1}{p}}E\left[\left(|S_n|-\varepsilon n^{\frac{1}{p}}\right)^{+}\right]<\infty\ \mbox{for all}\ \varepsilon>0,
\end{eqnarray}
where $x^+=\max\{0,x\}$.

Chow's result has been generalized in various directions. Wang and Su (2004), Wang et al. (2005), Chen (2006),  Guo and Xu (2006), Rosalsky et al. (2006), Ye and Zhu (2007), and Qiu et al. (2014) studied complete moment convergence for sums of Banach space valued random elements. Li and Zhang (20004), Chen et al. (2007), Kim et al. (2008), and Zhou (2010) considered complete moment convergence for moving average processes. Jiang and Zhang (2006), Li (2006), Liu and Lin (2006), Ye et al. (2007), Fu and Zhang (2008), Zhao and Tao (2008), and Chen and Zhang (2010) studied precise asymptotics for complete moment convergence. Wang and Zhao (2006), Liang et al. (2010), and Guo (2013) considered  complete moment convergence for negatively associated random variables. Qiu and Chen (2014) studied complete moment convergence for i.i.d. random variables, and extended two results in Gut and Stadtm\"{u}ller (2011) to complete moment convergence.

\begin{exa}
Let $\{X_n,n\geq 1\}$ be a sequence of  random variables with $\sup_{i\geq 1}E[|X_i|]\leq C$ for some positive constant $C$. Then for any $\alpha>1$, we have $\frac{S_n}{n^2(\ln n)^{\alpha}}\stackrel{s\mbox{-}L^1}{\longrightarrow} 0$.
In fact,
\begin{eqnarray*}
\sum_{n=1}^{\infty}E\left[\left|\frac{S_n}{n^2(\ln n)^{\alpha}}\right|\right]\leq \sum_{n=1}^{\infty}\frac{1}{n^2(\ln n)^{\alpha}}\sum_{k=1}^nE[|X_i|]\leq C\sum_{n=1}^{\infty}\frac{1}{n(\ln n)^{\alpha}}<\infty.
\end{eqnarray*}
\end{exa}

\begin{exa}
Let $\{X_n,n\geq 1\}$ be a sequence of pairwise uncorrelated random variables with $\sup_{i\geq 1}Var(X_i)\leq C$ for some positive constant $C$, where $Var(X_i)$ stands for the variance of $X_i$. Then for any $\alpha>1$, we have $\frac{S_n-E[S_n]}{n^{3/2}(\ln n)^{\alpha}}\stackrel{s^*\mbox{-}L^2}{\longrightarrow} 0$.
In fact,
\begin{eqnarray*}
\sum_{n=1}^{\infty}\left\|\frac{S_n-E[S_n]}{n^{3/2}(\ln n)^{\alpha}}\right\|_2&=&\sum_{n=1}^{\infty}\frac{1}{n^{3/2}(\ln n)^{\alpha}}(Var(S_n))^{1/2}\\
&=&\sum_{n=1}^{\infty}\frac{1}{n^{3/2}(\ln n)^{\alpha}}\left(\sum_{k=1}^nVar(X_i)\right)^{1/2}\\
&\leq& \sqrt{C}\sum_{n=1}^{\infty}\frac{1}{n(\ln n)^{\alpha}}<\infty.
\end{eqnarray*}
\end{exa}

In Section 3, we consider ``s-$L^p$ convergence " versions and ``s$^*$-$L^p$ convergence " versions of the Toeplitz lemma, the Ces\`{a}ro mean convergence theorem and the Kronecker lemma.
Four counterexamples will be  given  to show that
they can fail in general. Some sufficient conditions for the Ces\`{a}ro mean convergence theorem under these two complete moment convergences will be presented.

\section{Complete convergence}

\subsection{Counterexamples}
In this subsection, we will construct two counterexamples to show that ``complete convergence " versions of the Toeplitz lemma, the Ces\`{a}ro mean convergence theorem and the Kronecker lemma can fail in general.

The next example shows that complete convergence version of the Ces\`{a}ro mean convergence theorem fails.

\begin{exa}\label{exa2.1}
 Suppose that $\{X_n,n\geq 1\}$ is a sequence of independent random variables such that $P(X_n=n)=\frac{1}{n^2},P(X_n=0)=1-\frac{1}{n^2}$. For any $\varepsilon>0$, we have
$$
\sum_{n=1}^{\infty}P(|X_n-0|\geq \varepsilon)= \sum_{n=1}^{\infty}P(X_n=n)=\sum_{n=1}^{\infty}\frac{1}{n^2}<\infty,
$$
i.e. $X_n\to 0$ c.c. Let $\bar{X}_n=\frac{1}{n}\sum_{k=1}^nX_k,n\geq 1$. In the following, we will show that  $\bar{X}_n\nrightarrow 0$ c.c.

 Let $n=2k,k\geq 2$ and define $k$ sets $A_1,\cdots,A_k$ as follows:
\begin{eqnarray*}
&&A_1:=\{X_{2k}=2k\},\\
&&A_2:=\{X_{2k}=0,X_{2k-1}=2k-1\},\\
&&\ \ \ \ \ \cdots\\
&&A_k:=\{X_{2k}=0,\cdots,X_{k+2}=0,X_{k+1}=k+1\}.
\end{eqnarray*}
Then we have $\bigcup_{i=1}^kA_i\subset\{\bar{X}_n\geq \frac{1}{2}\}$, and thus
\begin{eqnarray*}
&&P\left(\bar{X}_n\geq \frac{1}{2}\right)\geq \sum_{i=1}^kP(A_i)\\
&&=\frac{1}{(2k)^2}+\left(1-\frac{1}{(2k)^2}\right)\frac{1}{(2k-1)^2}+\cdots+
\prod_{j=2}^k\left(1-\frac{1}{(k+j)^2}\right)\frac{1}{(k+1)^2}\\
&&\geq \prod_{j=2}^k\left(1-\frac{1}{(k+j)^2}\right)\sum_{i=k+1}^{2k}\frac{1}{i^2}.
\end{eqnarray*}
Denote $I_k= \prod_{j=2}^k\left(1-\frac{1}{(k+j)^2}\right)$. Then
\begin{eqnarray*}\label{exa2.1-a}
I_k&=&\frac{(2k+1)(2k-1)}{(2k)^2}\frac{2k(2k-2)}{(2k-1)^2}\cdots\frac{(k+4)(k+2)}{(k+3)^2}\frac{(k+3)(k+1)}{(k+2)^2}\\
&=&\frac{(2k+1)(k+1)}{2k(k+2)}\to 1\ \mbox{as}\ k\to\infty.
\end{eqnarray*}
Thus there exists a large number $K$ such that for any $k\geq K$, we have $I_k\geq \frac{1}{2}$. So, for any $n=2k\geq 2K$, we have
\begin{eqnarray*}
P\left(\bar{X}_n\geq \frac{1}{2}\right)\geq I_k\sum_{i=k+1}^{2k}\frac{1}{i^2}\geq\frac{1}{2}\sum_{i=k+1}^{2k}\frac{1}{(2k)^2}=\frac{1}{8k}.
\end{eqnarray*}
It follows that
$$
\sum_{n=1}^{\infty}P\left(\bar{X}_n\geq \frac{1}{2}\right)\geq \sum_{k=K}^{\infty}\frac{1}{8k}=\infty.
$$
Hence $\bar{X}_n\nrightarrow 0$ c.c.
\end{exa}

\begin{rem}
The above example also  shows the failure of the Toeplitz lemma when the mode of convergence is ``complete convergence", taking
$a_{nk}=1/n,1\leq k\leq k_n=n,n\geq 1$.
\end{rem}

The next example shows that complete convergence version of the Kronecker lemma also fails. The basic  idea comes from  Linero  and Rosalsky \cite[Example 2.3]{LR13}.

\begin{exa}
Let $\{Y_n,n\geq 1\}$ be a sequence of independent random variables such that
$P(Y_n=16^{n-1})=\frac{1}{n^2},P(Y_n=0)=1-\frac{1}{n^2}$. Denote $X_{2n-1}=Y_n,X_{2n}=-2Y_n,n\geq 1$. Then for any $n\geq 1$, we have
\begin{eqnarray}\label{Exa2.2-a}
\frac{X_{2n-1}}{2^{2n-1}}+\frac{X_{2n}}{2^{2n}}=0.
\end{eqnarray}
By the above definitions, for any $\varepsilon>0$,  we have
$$
\sum_{n=1}^{\infty}P(|X_n-0|\geq \varepsilon)\leq 2\sum_{n=1}^{\infty}\frac{1}{n^2}<\infty,
$$
 and thus $X_n\to 0$ c.c.  By (\ref{Exa2.2-a}), we know that  $\sum_{k=1}^n\frac{X_k}{2^k}=\frac{X_n}{2^n}I(n\ \mbox{is odd})$.
  Hence $\sum_{k=1}^n\frac{X_k}{2^k}\to 0$ c.c.

 In the following, we will show that $\frac{1}{2^n}\sum_{k=1}^nX_k\nrightarrow 0$ c.c. It's enough to show one of its subsequence
 \begin{eqnarray}\label{Exa2.2-b}
 \frac{1}{2^{4n}}\sum_{k=1}^{4n}X_k\nrightarrow 0\ c.c.
 \end{eqnarray}
 For any odd integer $k$,
 $$
 X_k+X_{k+1}=X_k-2X_k=-X_k.
 $$
 Thus, for any $n\geq 1$,
 $$
 \frac{1}{2^{2n}}\sum_{k=1}^{2n}X_k=-\frac{1}{2^{2n}}\sum_{k=1}^nX_{2k-1}.
 $$
 And so (\ref{Exa2.2-b}) can be expressed to be
 \begin{eqnarray}\label{Exa2.2-c}
 \frac{1}{16^n}\sum_{k=1}^{2n}X_{2k-1}\nrightarrow 0\ c.c.
 \end{eqnarray}
For  $k=n+1,\cdots,2n$,  we have
\begin{eqnarray}\label{Exa2.2-d}
P(X_{2k-1}=16^{k-1})=P(Y_k=16^{k-1})=\frac{1}{k^2},\ P(X_{2k-1}=0)=P(Y_k=0)=1-\frac{1}{k^2}.
\end{eqnarray}
 Define $n$ sets $A_1,\cdots,A_n$ as follows:
\begin{eqnarray*}
&&A_1:=\{X_{2(2n)-1}=16^{2n-1}\},\nonumber\\
&&A_2:=\{X_{2(2n)-1}=0,X_{2(2n-1)-1}=16^{(2n-1)-1}\},\nonumber\\
&&\ \ \ \ \cdots\nonumber\\
&&A_n:=\{X_{2(2n)-1}=0,\cdots,X_{2(n+2)-1}=0,X_{2(n+1)-1}=16^n\}.
\end{eqnarray*}
Then $\bigcup_{k=1}^nA_k\subset\{|\frac{1}{16^n}\sum_{k=1}^{2n}X_{2k-1}-0|\geq 1\}$, and thus
\begin{eqnarray}\label{Exa2.2-d}
&&P\left\{\left|\frac{1}{16^n}\sum_{k=1}^{2n}X_{2k-1}-0\right|\geq 1\right\}\geq \sum_{k=1}^nP(A_k)\nonumber\\
&&=\frac{1}{(2n)^2}+\left(1-\frac{1}{(2n)^2}\right)\frac{1}{(2n-1)^2}+\cdots+
\prod_{j=2}^n\left(1-\frac{1}{(n+j)^2}\right)\cdot\frac{1}{(n+1)^2}\nonumber\\
&&\geq \prod_{j=2}^n\left(1-\frac{1}{(n+j)^2}\right)\cdot\sum_{k=n+1}^{2n}\frac{1}{k^2}.
\end{eqnarray}
By (\ref{Exa2.2-d}) and following the deduction in Example 2.1, we can obtain that
$$
\sum_{n=1}^{\infty}P\left\{\left|\frac{1}{16^n}\sum_{k=1}^{2n}X_{2k-1}-0\right|\geq 1\right\}=\infty,
$$
i.e. (\ref{Exa2.2-c}) holds.
\end{exa}

\subsection{Sufficient conditions}

\begin{pro}\label{pro2.1}
Let $\{X_1,X_2,\cdots\}$ be pairwise uncorrelated random variables satisfying
\begin{eqnarray}\label{pro2.1-a}
\sum_{n=1}^{\infty}\frac{Var(X_n)}{n^{\alpha}}<\infty,
\end{eqnarray}
where $\alpha>0$, then for any $\varepsilon>0$,
\begin{eqnarray}\label{pro2.1-b}
\sum_{n=1}^{\infty}n^{1-\alpha}P\left\{\left|\frac{S_n-E(S_n)}{n}\right|\geq \varepsilon\right\}<\infty.
\end{eqnarray}
If  $\{X_1,X_2,\cdots\}$ is a sequence of independent random variables satisfying (\ref{pro2.1-a}), then for any $\varepsilon>0$,
\begin{eqnarray}\label{pro2.1-c}
\sum_{n=1}^{\infty}n^{1-\alpha}P\left\{\max_{1\leq k\leq n}|S_k-E(S_k)|\geq n\varepsilon\right\}<\infty.
\end{eqnarray}
\end{pro}
{\bf Proof.} For any $\varepsilon>0$, by  Chebyshev's inequality and (\ref{pro2.1-a}), we get
\begin{eqnarray}\label{pro2.1-d}
&&\sum_{n=1}^{\infty}n^{1-\alpha}P\left\{\left|\frac{S_n-E(S_n)}{n}\right|\geq \varepsilon\right\}\leq\sum_{n=1}^{\infty}n^{1-\alpha}\frac{Var(S_n)}{(n\varepsilon)^2}\nonumber\\
&&=\frac{1}{\varepsilon^2}\sum_{n=1}^{\infty}\frac{1}{n^{1+\alpha}}\sum_{k=1}^{n}Var(X_k)\nonumber\\
&&=\frac{1}{\varepsilon^2}\sum_{k=1}^{\infty}Var(X_k)\sum_{n=k}^{\infty}\frac{1}{n^{1+\alpha}}\nonumber\\
&&\leq\frac{M}{\varepsilon^2}\sum_{k=1}^{\infty}\frac{Var(X_k)}{k^{\alpha}}<\infty,
\end{eqnarray}
where $M$ is a positive constant. Hence (\ref{pro2.1-b}) holds. By Kolmogorov's inequality and the deduction of (\ref{pro2.1-d}), we get (\ref{pro2.1-c}).\hfill\fbox

\begin{rem}\label{rem2.2}
(i) Letting $\alpha=1$ in the above proposition, we get that if
$
\sum_{n=1}^{\infty}\frac{Var(X_n)}{n}<\infty$, then $\frac{S_n-E(S_n)}{n}\to 0$ c.c.

(ii) By Kolmogorov's strong law of large numbers, we know that if
$\{X_1,X_2,\cdots\}$ are independent random variables satisfying
$
\sum_{n=1}^{\infty}\frac{Var(X_n)}{n^2}<\infty,
$
then $\frac{S_n-E(S_n)}{n}\to 0$ a.s. By  Proposition \ref{pro2.1} and Baum and Katz \cite[Proposition 1(b)]{BK65}, we know
 that if $\{X_1,X_2,\cdots\}$ is a sequence of  pairwise uncorrelated random variables satisfying
$
\sum_{n=1}^{\infty}\frac{Var(X_n)}{n^2}<\infty,
$ and $|X_i|<i,\forall i\in \mathbb{N}$, then $\frac{S_n-E(S_n)}{n}\to 0$ a.s.
\end{rem}

\begin{cor}\label{cor2.3}
Let $\{X_1,X_2,\cdots\}$ be pairwise uncorrelated random variables satisfying \linebreak
$
\sum_{n=1}^{\infty}\frac{Var(X_n)}{n}<\infty,
$
and  $E(X_n)\to 0$. Then  $\frac{S_n}{n}\to 0$ c.c.
\end{cor}
{\bf Proof.} In this case, $\frac{E(S_n)}{n}=\frac{\sum_{k=1}^{n}E(X_k)}{n}\to 0$.
Then the result follows from Proposition \ref{pro2.1}.\hfill\fbox

\begin{rem}
By the above corollary, we know that if $\{X_1,X_2,\cdots\}$ is a sequence of  pairwise uncorrelated random variables satisfying that $X_n\to 0$ c.c. (or $X_n$ converges to 0 in probability),
$
\sum_{n=1}^{\infty}\frac{Var(X_n)}{n}<\infty,
$
and there exists an integrable random variable $X$ such that $|X_n|\leq X\ a.s.,\ \forall n\geq 1$. Then by the dominated convergence theorem, we have that $E(X_n)\to 0$ and thus by
the above corollary, the Ces\`{a}ro  sum $\frac{S_n}{n}$ of $\{X_n,n\geq 1\}$ satisfies
$\frac{S_n}{n}\to 0$ c.c.

\end{rem}

\section{Complete moment convergence}
\subsection{Counterexamples}

In this subsection, we will construct four counterexamples to show that s-$L^1$ convergence versions and s$^*$-$L^2$ convergence versions  of the Toeplitz lemma, the Ces\`{a}ro mean convergence theorem and the Kronecker lemma can fail in general.

The next example shows that s-$L^1$ convergence versions of the Ces\`{a}ro mean convergence theorem and the Toeplitz lemma fail.

\begin{exa}\label{exa3.1}
 Let $\{X_n,n\geq 1\}$ be a sequence of random variables such that $P(X_n=n)=\frac{1}{n^3},P(X_n=0)=1-\frac{1}{n^3}$.
Then we have $E[|X_n|]=\frac{1}{n^2}$ and thus
\begin{eqnarray*}
\sum_{n=1}^{\infty}E[|X_n|]=\sum_{n=1}^{\infty}\frac{1}{n^2}<\infty,
\end{eqnarray*}
i.e. $X_n\stackrel{s\mbox{-}L^1}{\longrightarrow}0$.

Let $\bar{X}_n=\frac{1}{n}\sum_{k=1}^nX_k,n\geq 1$. Then
\begin{eqnarray*}
E[|\bar{X}_n|]=\frac{1}{n}\sum_{k=1}^nE[|X_k|]=\frac{1}{n}\sum_{k=1}^n\frac{1}{k^2}.
\end{eqnarray*}
Since $\sum_{k=1}^{\infty}\frac{1}{k^2}=\frac{\pi^2}{6}$, there exists a large $N$ such that $\forall n\geq N,\sum_{k=1}^n\frac{1}{k^2}\geq \frac{\pi^2}{12}$. Hence
$$
\sum_{n=1}^{\infty}E[|\bar{X}_n|]\geq \sum_{n=N}^{\infty}\frac{1}{n}\cdot \frac{\pi^2}{12}=\infty,
$$
and so it  doesn't hold that $\bar{X}_n\stackrel{s\mbox{-}L^1}{\longrightarrow} 0$.
\end{exa}

The next example shows that s-$L^1$ convergence version of the Kronecker lemma fails. The basic  idea comes from  Linero  and Rosalsky \cite[Example 2.3]{LR13}.

\begin{exa}\label{exa3.2}
Let $\{Y_n,n\geq 1\}$ be a sequence of independent random variables such that
$P(Y_n=n)=\frac{1}{n^2(\ln n)^{1+\alpha}},P(Y_n=0)=1-\frac{1}{n^2(\ln n)^{1+\alpha}}$,  $\alpha>0$. Denote $X_{2n-1}=(2n-1)Y_n,X_{2n}=-2nY_n,n\geq 1$. Then for any $n\geq 1$, we have
\begin{eqnarray}\label{Exa3.2-a}
\frac{X_{2n-1}}{2n-1}+\frac{X_{2n}}{2n}=0.
\end{eqnarray}
 By (\ref{Exa3.2-a}), we know that  $\sum_{k=1}^n\frac{X_k}{k}=\frac{X_n}{n}I(n\ \mbox{is odd})$.
 If $n=2k-1$, then we have
 $$
 E[|X_n/n|]=E[|Y_k|]=\frac{1}{k(\ln k)^{1+\alpha}},
 $$
 which implies that
 $$
 \sum_{n=1}^{\infty}E\left[\left|\sum_{k=1}^n\frac{X_k}{k}\right|\right]= \sum_{k=1}^{\infty}\frac{1}{k(\ln k)^{1+\alpha}}<\infty,
 $$
 i.e. $\sum_{k=1}^n\frac{X_k}{k}\stackrel{s\mbox{-}L^1}{\longrightarrow} 0$.

 In the following, we will show that $\frac{1}{n}\sum_{k=1}^nX_k\stackrel{s\mbox{-}L^1}{\nrightarrow} 0$. It's enough to show one of its subsequence
 \begin{eqnarray}\label{Exa3.2-b}
 \frac{1}{2n}\sum_{k=1}^{2n}X_k\stackrel{s\mbox{-}L^1}{\nrightarrow} 0.
 \end{eqnarray}
 For any  integer $k$, we have $ X_{2k-1}+X_{2k}=-Y_k.$
 Thus (\ref{Exa3.2-b}) can be expressed to be
 \begin{eqnarray}\label{Exa3.2-c}
 \frac{1}{2n}\sum_{k=1}^{n}Y_k\stackrel{s\mbox{-}L^1}{\nrightarrow} 0.
 \end{eqnarray}
 By the Fubini theorem, we have
 \begin{eqnarray*}
 \sum_{n=1}^{\infty}E\left[\left|\frac{1}{2n}\sum_{k=1}^{n}Y_k\right|\right]
 &= &\sum_{n=1}^{\infty}\frac{1}{2n}\sum_{k=1}^{n}E\left[Y_k\right]
 =\sum_{n=1}^{\infty}\frac{1}{2n}\sum_{k=1}^{n}\frac{1}{k(\ln k)^{1+\alpha}}\\
 &=&\sum_{k=1}^{\infty}\frac{1}{k(\ln k)^{1+\alpha}}\sum_{n=k}^{\infty}\frac{1}{2n}=\infty.
 \end{eqnarray*}
 Hence (\ref{Exa3.2-c}) holds.

\end{exa}


The next example shows that s$^*$-$L^2$ convergence versions of the Ces\`{a}ro mean convergence theorem and the Toeplitz lemma fail.

\begin{exa}\label{exa3.3}
Let $\{X_n,n\geq 1\}$ be a sequence of random variables such that $P(X_n=\sqrt{n})=\frac{1}{n^5},P(X_n=0)=1-\frac{1}{n^5}$.
Then we have $E[|X_n|^2]=\frac{1}{n^4}$ and thus
\begin{eqnarray*}
\sum_{n=1}^n\|X_n\|_2=\sum_{n=1}^{\infty}\left(\frac{1}{n^4}\right)^{\frac{1}{2}}=\sum_{n=1}^{\infty}\frac{1}{n^2}<\infty,
\end{eqnarray*}
i.e. $X_n\stackrel{s^*\mbox{-}L^2}{\longrightarrow} 0$.

Let $\bar{X}_n=\frac{1}{n}\sum_{k=1}^nX_k,n\geq 1$. Then
\begin{eqnarray*}
E[|\bar{X}_n|^2]&=&\frac{1}{n^2}\left(\sum_{k=1}^nE[|X_k|^2]+2\sum_{1\leq i<j\leq n}E[X_iX_j]\right)\\
&\geq&\frac{1}{n^2}\sum_{k=1}^nE[|X_k|^2]=\frac{1}{n^2}\sum_{k=1}^n\frac{1}{k^4}.
\end{eqnarray*}
Denote $c=\sum_{k=1}^{\infty}\frac{1}{k^4}$. Then $c$ is a positive constant and there exists a large $N$ such that
$\forall n\geq N,\sum_{k=1}^n\frac{1}{k^4}\geq \frac{c}{2}$. It follows that
\begin{eqnarray*}
\sum_{n=1}^{\infty}\|\bar{X}_n\|_2\geq \sum_{n=N}^{\infty}\left(\frac{1}{n^2}\cdot \frac{c}{2}\right)^{\frac{1}{2}}=\sqrt{\frac{c}{2}}\sum_{n=N}^{\infty}\frac{1}{n}=\infty.
\end{eqnarray*}
Hence it doesn't hold that $\bar{X}_n\longrightarrow^{\!\!\!\!\!\!\!\!\!\!\! \!\!\!s^*\mbox{-}L^2} 0$.
\end{exa}

Following Examples \ref{exa3.2} and \ref{exa3.3}, we construct the following example, which shows that  s$^*$-$L^2$ convergence version of the Kronecker lemma fails.

\begin{exa}\label{exa3.4}
Let $\{Y_n,n\geq 1\}$ be a sequence of independent random variables such that
$P(Y_n=\sqrt{n})=\frac{1}{n^5},P(Y_n=0)=1-\frac{1}{n^5}$. Denote $X_{2n-1}=(2n-1)Y_n,X_{2n}=-2nY_n,n\geq 1$. Then for any $n\geq 1$, we have
\begin{eqnarray}\label{Exa3.4-a}
\frac{X_{2n-1}}{2n-1}+\frac{X_{2n}}{2n}=0.
\end{eqnarray}
 By (\ref{Exa3.4-a}), we know that  $\sum_{k=1}^n\frac{X_k}{k}=\frac{X_n}{n}I(n\ \mbox{is odd})$.
 If $n=2k-1$, then we have
 $$
 \|X_n/n\|_2=\|Y_k\|_2=\frac{1}{k^2}.
 $$
 Hence
 $$
 \sum_{n=1}^{\infty}\left\|\sum_{k=1}^n\frac{X_k}{k}\right\|_2=\sum_{k=1}^{\infty}\|Y_k\|_2=\sum_{k=1}^{\infty}\frac{1}{k^2}<\infty,
 $$
i.e. $\sum_{k=1}^n\frac{X_k}{k}\stackrel{s^*\mbox{-}L^2}{\longrightarrow} 0$.

 In the following, we will show that $\frac{1}{n}\sum_{k=1}^nX_k\stackrel{s^*\mbox{-}L^2}{\nrightarrow} 0$. It's enough to show one of its subsequence
 \begin{eqnarray}\label{Exa3.4-b}
 \frac{1}{2n}\sum_{k=1}^{2n}X_k\stackrel{s^*\mbox{-}L^2}{\nrightarrow} 0.
 \end{eqnarray}
 For any  integer $k$, we have $ X_{2k-1}+X_{2k}=-Y_k.$
 Thus (\ref{Exa3.4-b}) can be expressed to be
 \begin{eqnarray}\label{Exa3.4-c}
 \frac{1}{2n}\sum_{k=1}^{n}Y_k\stackrel{s^*\mbox{-}L^2}{\nrightarrow} 0.
 \end{eqnarray}
 Denote $c=\sum_{k=1}^{\infty}\frac{1}{k^4}$. Then $0<c<\infty$, and there exists $N$ such that for any $n\geq N$, we have
 $\sum_{k=1}^{n}\frac{1}{k^4}\geq \frac{c}{2}$. Hence we have
\begin{eqnarray*}
\sum_{n=1}^{\infty}\left\|\frac{1}{2n}\sum_{k=1}^{n}Y_k\right\|_2&=&\sum_{n=1}^{\infty}\frac{1}{2n}\left(\sum_{k=1}^{n}E[Y_k^2]+2\sum_{1\leq i<j\leq n}E[Y_iY_j]\right)^{1/2}\\
&\geq&\sum_{n=1}^{\infty}\frac{1}{2n}\left(\sum_{k=1}^{n}E[Y_k^2]\right)^{1/2}=\sum_{n=1}^{\infty}\frac{1}{2n}\left(\sum_{k=1}^{n}\frac{1}{k^4}\right)^{1/2}\\
&\geq&\sum_{n=N}^{\infty}\frac{1}{2n}\left(\sum_{k=1}^{n}\frac{1}{k^4}\right)^{1/2}\\
&\geq&\sum_{n=N}^{\infty}\frac{1}{2n}\sqrt{\frac{c}{2}}=\infty.
\end{eqnarray*}
Hence (\ref{Exa3.4-c}) holds.
\end{exa}

\subsection{Sufficient conditions}

By  Example 3.1, we know that, if $\sum_{n=1}^{\infty}E[|X_n|^p]<\infty$, then we
don't have $\sum_{n=1}^{\infty}E[|S_n/n|^p]<\infty$ necessarily. In general, we have the following result.

\begin{pro}\label{pro3.1}
Suppose that $1\leq p< \infty$ and $\sum_{n=1}^{\infty}E[|X_n|^p]<\infty$, then  $\forall \varepsilon>0$, we have
\begin{eqnarray*}
\sum_{n=1}^{\infty}\frac{1}{(\ln n)^{1+\varepsilon}}E\left[\left|S_n/n\right|^p\right]<\infty.
\end{eqnarray*}
\end{pro}
{\bf Proof.} By the convexity of  the function $f(x)=|x|^p$, we have
\begin{eqnarray*}
\sum_{n=1}^{\infty}\frac{1}{(\ln n)^{1+\varepsilon}}E\left[\left|S_n/n\right|^p\right]&\leq& \sum_{n=1}^{\infty}\frac{1}{n(\ln n)^{1+\varepsilon}}\left(\sum_{k=1}^nE[|X_k|^p]\right)\\
&\leq &\left(\sum_{k=1}^{\infty}E[|X_k|^p]\right)\sum_{n=1}^{\infty}\frac{1}{n(\ln n)^{1+\varepsilon}}<\infty.
\end{eqnarray*}
 \hfill\fbox

 \begin{pro}\label{pro3.2}
Let  $\{X_1,X_2,\cdots\}$ be pairwise uncorrelated random variables satisfying \linebreak
$\sum_{n=1}^{\infty}Var(X_n)<\infty,$
 then for any $1<q\leq 2$, we have
 \begin{eqnarray*}
 \sum_{n=1}^{\infty}E\left[\left|\frac{S_n-E(S_n)}{n}\right|^q\right]<\infty,
 \end{eqnarray*}
 in particular, $\frac{S_n-E(S_n)}{n}\to 0$ c.c.
 \end{pro}
{\bf Proof.} By the assumptions, we have
\begin{eqnarray*}
\sum_{n=1}^{\infty}E\left[\left|\frac{S_n-E(S_n)}{n}\right|^q\right]&=&\sum_{n=1}^{\infty}\left(\left\|\frac{S_n-E(S_n)}{n}\right\|_q\right)^q\\
&\leq&\sum_{n=1}^{\infty}\left(\left\|\frac{S_n-E(S_n)}{n}\right\|_2\right)^q\\
&=&\sum_{n=1}^{\infty}\frac{1}{n^q}\left(\sum_{i=1}^{n}Var(X_i)\right)^{q/2}\\
&\leq&\sum_{n=1}^{\infty}\frac{1}{n^q}\left(\sum_{i=1}^{\infty}Var(X_i)\right)^{q/2}<\infty.
\end{eqnarray*}
\hfill\fbox

By Example \ref{exa3.3}, we know that, if $\sum_{n=1}^{\infty}\|X_n\|_p<\infty$, then we
don't have $\sum_{n=}^{\infty}\|\frac{S_n}{n}\|_p<\infty$ necessarily. In general, we have the following two propositions.

\begin{pro}\label{pro3.3}
Suppose that $1\leq p< \infty$ and $\sum_{n=1}^{\infty}\|X_n\|_p<\infty$, then  $\forall \varepsilon>0$, we have
\begin{eqnarray}\label{pro3.3-a}
\sum_{n=1}^{\infty}\frac{1}{(\ln n)^{1+\varepsilon}}\left\|S_n/n\right\|_p<\infty.
\end{eqnarray}
\end{pro}
{\bf Proof.} By  Minkowski's inequality and the definition of the norm $\|\cdot\|_p$, we have that
$$
\left\|S_n/n\right\|_p\leq \frac{1}{n}(\sum_{k=1}^{n}\|X_k\|_p).
$$
Then we can prove (\ref{pro3.3-a}) by following the proof of Proposition \ref{pro3.1}.\hfill\fbox

\begin{pro}\label{pro3.4}
Suppose that $1<p<\infty$ and
$
\sum_{n=1}^{\infty}\|X_n\|_p<\infty,
$
then for any $1<q\leq p$, we have
\begin{eqnarray*}
\sum_{n=1}^{\infty}E[|S_n/n|^q]<\infty,
\end{eqnarray*}
in particular, $S_n/n\to 0$ c.c.
\end{pro}
{\bf Proof.} By the fact that $\|\cdot\|_q\leq \|\cdot\|_p$, Minkowski's inequality and the
assumption, we have
\begin{eqnarray*}
\sum_{n=1}^{\infty}E[|S_n/n|^q]&=&\sum_{n=1}^{\infty}\left(\|S_n/n\|_q\right)^q\leq\sum_{n=1}^{\infty}\left(\|S_n/n\|_p\right)^q\\
&\leq&\sum_{n=1}^{\infty}\left(\frac{\sum_{k=1}^n\|X_k\|_p}{n}\right)^q\\
&=&\sum_{n=1}^{\infty}\frac{1}{n^q}\left(\sum_{k=1}^n\|X_k\|_p\right)^q\\
&\leq&\left(\sum_{k=1}^{\infty}\|X_k\|_p\right)^q\sum_{n=1}^{\infty}\frac{1}{n^q}<\infty.
\end{eqnarray*}
\hfill\fbox

\begin{pro}\label{pro3.5}
Suppose that $\sum_{n=1}^{\infty}\|X_n\|_{\infty}<\infty$. Then\\
(i) for any $\varepsilon>0$, we have
\begin{eqnarray}\label{pro3.4-a}
\sum_{n=1}^{\infty}\frac{1}{(\ln n)^{1+\varepsilon}}\|S_n/n\|_{\infty}<\infty;
\end{eqnarray}
(ii) for any $1<q<\infty$, we have
\begin{eqnarray*}\label{pro3.4-b}
\sum_{n=1}^{\infty}E[|S_n/n|^q]<\infty,
\end{eqnarray*}
in particular, $S_n/n\to 0$ c.c.
\end{pro}
{\bf Proof.} (i) By the definition of the norm $\|\cdot\|_{\infty}$, we have that
$$
\|S_n/n\|_{\infty}\leq \frac{1}{n}\left(\sum_{k=1}^n\|X_k\|_{\infty}\right).
$$
Then we can prove (\ref{pro3.4-a}) by following the proof of Proposition \ref{pro3.1}.

(ii) It's  a direct consequence of Proposition \ref{pro3.4} by noting that for any $1<p<\infty$ and any
random variable $X$, $\|X\|_p\leq \|X\|_{\infty}$.\hfill\fbox

\bigskip

{ \noindent {\bf\large Acknowledgments} \vskip 0.1cm  \noindent  We acknowledge the helpful suggestions and comments of three anonymous referees,
which improved the presentation of this paper.  We are grateful to the support of
NNSFC (Grant No. 11371191) and Jiangsu Province Basic Research Program (Natural Science
Foundation) (Grant No. BK2012720).}


\begin{thebibliography}{1234}

\bibitem{BK65} Baum, L. E., Katz, M. (1965). Convergence rates  in the law of large numbers.  Trans. Amer. Math. Soc. 120:108-123.

\bibitem{Ch06} Chen, P. (2006). Complete moment convergence for sequences of independent random elements in Banach spaces.  Stoch. Anal. Appl. 24:999-1010.

\bibitem{CZ10}Chen, Y.-Y., Zhang, L.-X. (2010). Second moment convergence rates for uniform empirical processes.  J. Inequal. Appl. 2010:972324, 9 pages (doi: 10.1155/2010/972324).

\bibitem{CHV07} Chen, P., Hu. T.-C., Volodin, A. (2007). Limiting behavior of moving average processes under negative association.  Teor. Imovir. ta Matem. Statyst. 7:154-166.

\bibitem{CS14} Chen, P., Sung S. H. (2014). A Baum-Katz theorem for i.i.d. random variables with higher order moments. Statis. Probab. Lett. 94:63-68.

\bibitem{Ch88} Chow, Y. (1988). On the rate of moment convergence of sample sums and extremes. Bull. Inst. Math. Acad. Sin. 16:177-201.

\bibitem{Er49} Erd\"{o}s, P. (1949). On a theorem of Hsu and Robbins.  Ann. Math. Statist. 20:286-291.

\bibitem{Er50} Erd\"{o}s, P. (1950). Remark on my paper ``On a theorem of Hsu and Robbins".  Ann. Math. Statist. 21:138.

\bibitem{FZ08} Fu, K.-A., Zhang, L.-X. (2008). On the moment convergence rates of LIL in Hilbert space. Math. Comput. Modelling 47:153-167.

\bibitem{Guo13}Guo M. (2013). On complete moment convergence of weighted sums for arrays of row-wise negatively associated random variables. Stochastics An International Journal of Probability and Stochastic Processes DOI: 10.1080/17442508.2013.801971.

\bibitem{GX06} Guo, M.-L., Xu, J. (2006). Moment complete convergence for arrays of $B$-valued random elements (in Chinese).  Appl. Math. J. Chinese Univ. Ser. A 21:432-438.


\bibitem{Gu78} Gut, A. (1978). Marcinkiwicz laws and convergence rates in the law of large numbers for random variables with multidimensionial indices. Ann. Probab. 6:469-482.

\bibitem{Gu80} Gut, A. (1980). Convergence rates for probabilities of moderate deviations for sums of random variables with multidimensional indices.  Ann. Probab. 8:298-313.

\bibitem{GS11} Gut, A., Stadtm\"{u}ller, U. (2011). An intermediate Baum-Katz theorem.  Statis. Probab. Lett. 81:1486-1492.


\bibitem{HR47}Hsu, P.,  Robbins, H. (1947). Complete convergence and the law of large numbers.  Proc. Natl. Acad. Sci. USA 33:25-31.

\bibitem{JZ06} Jiang, Y., Zhang, L.-X. (2006). Precise asymptotics in complete moment convergence of i.i.d. random variables (in Chinese). Acta. Math. Sci. Ser. A Chin. Ed. 26:917-925.


\bibitem{KKC08}Kim T.-S., Ko M.-H., Choi Y.-K. (2008). Complete moment convergence of moving average processes with dependent innovations. J. Korean Math. Soc. 45:355-365.

\bibitem{La98} Lanzinger, H. (1998). A Baum-Katz theorem for random variables under exponential momment conditions. Statis. Probab. Lett. 39:89-95.

\bibitem{LRJW95} Li, D., Rao, M.B., Jiang, T., Wang X. (1995). Complete convergence and almost sure convergence of weighted sums of random variables.  J. Theoret. Probab. 8:49-76.

\bibitem{Li06}Li Y.-X. (2006). Precise asymptotics in complete moment convergence of moving-average processes. Statis. Probab. Lett. 76:1305-1315.

\bibitem{LZ04} Li, Y.-X., Zhang, L.-X. (2004). Complete moment convergence of moving-average processes under dependence assumptions.  Statis. Probab. Lett. 70:191-197.

\bibitem{LLR10} Liang, H. Y., Li, D. L., Rosalsky, A. (2010). Complete moment and integral convergence for sums of negatively associated random variables.  Acta Math. Sin. Engl. Ser. 26:419-432.

 \bibitem{LR13} Linero A., Rosalsky A. (2013).  On the Toeplitz lemma, convergence in probability, and mean convergence.   Stoch. Anal. Appl. 31:684-694.

\bibitem{LL06} Liu, W., Lin Z. (2006). Precise asymptotics for a new kind of complete moment convergence. Statis. Probab. Lett. 76:1787-1799.

\bibitem{L77} Lo\`{e}ve, M. (1977).  {\it Probability Theory I (4th ed.)}  New York: Springer-Verlag.

\bibitem{QC14}Qiu, D., Chen P. (2014). Complete moment convergence for i.i.d. random variables. Statis. Probab. Lett. 91:76-82.

\bibitem{QUV14}Qiu, D., Urmeneta, H., Volodini, A. (2014). Complete moment convergence for weighted sums of sequences of independent random elements in Banach spaces.  Collect. Math.  65:155-167.


\bibitem{RTV06} Rosalsky, A., Thanh, L. V., Volodin, A. (2006). On complete convergence in mean of normed sums of independent random elements in Banach spaces.  Stoch. Anal. Appl. 24:23-35.

\bibitem{Su07}Sung, S. H. (2007). Complete convergence for weighted sums of random variables. Statis. Probab. Lett. 77:303-311.

\bibitem{SV06}Sung, S.H., Volodin, A. (2006). On the rate of complete convergence for weighted sums of arrays of random elements. J. Korean Math. Soc. 43:815-828.

\bibitem{WS04} Wang, D., Su, C. (2004). Moment complete convergence for sequences of B-valued iid random elements (in Chinese).  Acta Math. Appl. Sin. 27:440-448.

\bibitem{WZW05} Wang, Y., Zhong, S. M., Wang, D. C. (2005). Moment complete convergence for sums of i.i.d. random elements in Banach space (in Chinese).  Dianzi Keji Daxue Xuebao 34:410-412.


\bibitem{WZ06}Wang, D., Zhao, W. (2006). Moment complete convergence for sums of a sequence of NA random variables (in Chinese).  Appl. Math. J. Chinese Univ. Ser. A 21:445-450.

\bibitem{YZ07}Ye, F., Zhu, D. J. (2007). Moment complete convergence for $B$-valued random elements (in Chinese).  J. Hefei Univ. Technol. Nat. Sci. 30:525-528.

\bibitem{YZP07} Ye, J., Zhang L.-X., Pang, T.-X. (2007). Precise rates in the law of logarithm for the moment convergence of i.i.d. random variables. J. Math. Anal. Appl. 327:695-714.

\bibitem{Z10}Zhou X. (2010). Complete moment convergence of moving average processes under $\varphi$-mixing assumptions. Statis. Probab. Lett. 80:285-292.

\bibitem{ZT08}Zhao Y., Tao J. (2008). Precise asymptotics in complete moment convergence for self-normalized sums. Comptut. Math. Appl. 56:1779-1786.

\end{thebibliography}
\end{document}